\def\draft{n}
\documentclass{amsart}
\usepackage[headings]{fullpage}
\usepackage{amssymb,epic,eepic,epsfig,amsbsy,amsmath}

\theoremstyle{plain}

\newtheorem{theorem}{Theorem}
\newtheorem{proposition}{Proposition}[section]

\theoremstyle{definition}
\newtheorem{definition}[proposition]{Definition}

\theoremstyle{remark}

\newtheorem{remark}[proposition]{Remark}

\def\printname#1{
    \if\draft y
        \smash{\makebox[0pt]{\hspace{-0.5in}
            \raisebox{8pt}{\tt\tiny #1}}}
    \fi
}

\newcommand{\psdraw}[2]
         {\begin{array}{c} \hspace{-1.3mm}
    \raisebox{-4pt}{\epsfig{figure=draws/#1.eps,width=#2}}
    \hspace{-1.9mm}\end{array}}

\newlength{\standardunitlength}
\catcode`\@=11
\long\def\@makecaption#1#2{%
     \vskip 10pt

\setbox\@tempboxa\hbox{
       \small\sf{\bfcaptionfont #1. }\ignorespaces #2}%
     \ifdim \wd\@tempboxa >\captionwidth {%
         \rightskip=\@captionmargin\leftskip=\@captionmargin
         \unhbox\@tempboxa\par}%
       \else
         \hbox to\hsize{\hfil\box\@tempboxa\hfil}%
     \fi}
\font\bfcaptionfont=cmssbx10 scaled \magstephalf
\newdimen\@captionmargin\@captionmargin=2\parindent
\newdimen\captionwidth\captionwidth=\hsize
\catcode`\@=12

\def\lbl#1{\label{#1}\printname{#1}}

\def\BN{\mathbb N}
\def\BZ{\mathbb Z}

\def\BQ{\mathbb Q}
\def\BR{\mathbb R}
\def\BC{\mathbb C}

\def\a{\alpha}

\def\La{\Lambda}
\def\l{\lambda}
\def\Ga{\Gamma}

\def\ga{\gamma}

\def\w{\omega}

\def\e{\epsilon}

\def\Ga{\Gamma}
\def\d{\delta}
\def\b{\beta}

\def\longto{\longrightarrow}

\def\w{\omega}

\def\calL{\mathcal{L}}

\def\calB{\mathcal{B}}

\def\calN{\mathcal{N}}
\def\Res{\mathrm{Res}}
\def\Rhat{\hat{R}}

\def\pp{\mathfrak{t}}
\def\Lnp{L^{\mathrm{NP}}}
\def\Lp{L^{\mathrm{P}}}

\begin{document}


\title[Resurgence of the Euler-MacLaurin summation formula]{
Resurgence of the Euler-MacLaurin summation formula}

\author{Ovidiu Costin}
\address{Department of Mathematics \\
         Ohio State University \\ 
         231 W 18th Avenue \\
         Columbus, OH 43210, USA \\ \newline
         {\tt http://www.math.ohio-state.edu/$\sim$costin } }
\email{costin@math.ohio-state.edu}
\author{Stavros Garoufalidis}
\address{School of Mathematics \\
         Georgia Institute of Technology \\
         Atlanta, GA 30332-0160, USA \\ \newline
         {\tt http://www.math.gatech.edu/$\sim$stavros } }
\email{stavros@math.gatech.edu}

\thanks{O.C. was supported in part by NSF grants DMS-0406193 and 
DMS-0600369 and S.G was supported in part by NSF grant DMS-0505445. \\
\newline
1991 {\em Mathematics Classification.} Primary 57N10. Secondary 57M25.
}

\date{ July 9, 2007 } 


\begin{abstract}
The Euler-MacLaurin summation formula relates a sum of a function to
a corresponding integral, with a remainder term. The remainder term has 
an asymptotic expansion, and for a typical analytic function, it
is a divergent (Gevrey-1) series. Under some decay assumptions of the function 
in a half-plane (resp. in the vertical strip containing 
the summation interval),
Hardy (resp. Abel-Plana) prove that the asymptotic expansion is a
Borel summable series, and give an exact Euler-MacLaurin summation formula.

Using a mild resurgence hypothesis for the function to be summed, we
give a Borel summable transseries expression for the remainder term, 
as well as a Laplace integral formula, with an explicit 
integrand which is a resurgent function itself. In particular, our summation
formula allows for resurgent functions with singularities in the vertical
strip containing the summation interval. 

Finally, we give two applications of our results. One concerns the
construction of solutions of linear difference equations with a small 
parameter. And another concerns resurgence of 1-dimensional sums
of quantum factorials, that are associated to knotted 3-dimensional objects.
\end{abstract}

\maketitle

\tableofcontents


\section{Introduction}
\lbl{sec.intro}

\subsection{The Euler-MacLaurin summation formula}
\lbl{sub.maclaurin}

The {\em Euler-MacLaurin summation formula} relates summation to integration
in the following way (see \cite[Sec.8]{O}):

\begin{equation}
\lbl{eq.em2}
\sum_{k=1}^N f\left(\frac{k}{N}\right)=
N \int_0^1 f(s) ds + \frac{1}{2}(f(1)-f(0)) + R(f,N)
\end{equation}
where the remainder $R(f,N)$ has an asymptotic expansion
\begin{equation}
\lbl{eq.RfN2}
R(f,N) \sim \Rhat(f,N)
\end{equation}
in the sense of Poincar\'e, where
\begin{equation}
\lbl{eq.RfN3}
\Rhat(f,N)= \sum_{n=1}^\infty \frac{B_{2n}}{(2n)!} \left( f^{(2n-1)}(1)-
f^{(2n-1)}(0) \right) \frac{1}{N^{2n-1}} \in \BC[[N^{-1}]].
\end{equation}
and $B_m$ are the {\em Bernoulli numbers} defined by the generating
series
\begin{equation}
\lbl{eq.bernoulli}
\frac{p}{e^p-1}=\sum_{n=0}^\infty \frac{B_n}{n!} p^n.
\end{equation}

Typically, the formal power series $\Rhat(f,N)$ is divergent and Gevrey-1,
due to the fact that the $n$-th derivative in \eqref{eq.RfN3} 
is not divided by an $n!$. In the present paper, we 
discuss an exact form of the Euler-MacLaurin summation
formula, under a resurgence hypothesis of the function $f(x)$; see
Proposition \ref{prop.2} below.

\subsection{Two applications of our exact Euler-MacLaurin summation formula}
\lbl{sub.applications}

Our exact form of the Euler-MacLaurin summation formula has two applications:
in {\em Quantum Topology} 
(where one sometimes needs to apply the Euler-MacLaurin 
summation formula to a resurgent function that has singularities in
the vertical strip which is perpendicular to the range of summation),
and in Borel summability (with respect to $\e)$ of difference equations
with a small $\e$-parameter. Let us discuss these applications.

Consider a triple $\pp=(a,b,\e)$ where $a,b\in \BZ$, $b > 0$, $\e=\pm 1$
and the expression 
 
\begin{equation}
\lbl{eq.Fabe}
I_{\pp}(q)=\sum_{n=0}^\infty q^{a \frac{n(n+1)}{2}} (q)^b_n \e^{ n}
\end{equation}
where $(q)_n$ is the {\em quantum factorial} defined by:
\begin{equation}
\lbl{eq.qfac}
(q)_n=\prod_{k=1}^n (1-q^k), \qquad (q)_0=1.
\end{equation}
Although $I_{\pp}(q)$ does not makes sense when $q$ is inside or outside the
unit disk, it does makes sense when 
\begin{itemize}
\item[(a)]
$q$ is a complex root of unity; in that case $I_{\pp}(q) \in \BC$.
\item[(b)]
$q=e^{1/x}$; in that case $I_{\pp}(q) \in \BQ[[1/x]]$.
\end{itemize}
Given $\pp$ as above, consider the power series:
\begin{eqnarray}
\lbl{eq.LNP}
\Lnp_{\pp}(p) &=& 1+\sum_{n=0}^\infty I_{\pp}(e^{\frac{2 \pi i}{n}}) p^n \\
\lbl{eq.LP}
\Lp_{\pp}(p) &=&  \calB( I_{\pp}(e^{1/x}))
\end{eqnarray}
where $\calB$ is the Borel transform defined below in Definition 
\ref{def.borelt}.
Our present results, together with \cite{CG3} and additional arguments, 
imply the following theorem, which will be presented in detail in 
a forthcoming publication \cite{CG2}.

\begin{theorem}
\lbl{thm.CG2}\cite{CG2}
For all $\pp$ as above, the power series $\Lnp_{\pp}(p)$ and $\Lp_{\pp}(p)$
are resurgent functions.
\end{theorem}
In particular, it follows that the generating series 
$\Lnp_{\pp}(p)$ of the Kashaev invariants of the two simplest knots, 
the $3_1$ and $4_1$ (corresponding to $\pp=(0,1,1)$ and $\pp=(2,-1,-1)$)
are resurgent functions.

Another application of our exact Euler-MacLaurin formula is to 
prove parametric resurgence (i.e., resurgence with respect
to $\e$) of a formal (WKB) solution to a linear difference equation with a 
small
parameter:
\begin{equation}
\lbl{eq.diffa}
y(x+\e,\e)=a(x,\e) y(x,\e)
\end{equation}
The formal solution of \eqref{eq.diffa} is of the form:
\begin{equation}
\lbl{eq.WKB2}
y(x,\e)=\exp\left(\frac{1}{\e} \sum_{k=0}^\infty F_k(x) \e^k \right).
\end{equation}
Under suitable hypothesis on $a(x,\e)$, Theorem \ref{thm.diff} below proves 
resurgence of the above series for
$x$ fixed, and constructs an actual solution to \eqref{eq.diffa} which
is asymptotic to the formal solution \eqref{eq.WKB2}. 

\subsection{Known forms of the Euler-MacLaurin summation formula}
\lbl{sub.knownforms}

Before we state our results, let us recall what is already
known. Suppose that $f$ satisfies the following assumption:

\begin{itemize}
\item[]
$f$ is analytic and satisfies the following bound:
\begin{equation}
\lbl{eq.hardyb}
f(x)=O(|x|^{-s})
\end{equation}
for some $0 < \d < 1$ and $s>0$, uniformly in the right-half plane 
$\Re(x) \geq \d$. 
\end{itemize}
For such functions $f$, 
Hardy proved in \cite[Sec.13.15]{Ha} that $\Rhat(f,N)$ is {\em
Borel summable}, and that the Borel sum agrees with the original sum.
In other words, the Borel tranform of $\Rhat(f,N)$ can be extended to the 
ray $[0,\infty)$, it is integrable of at most exponential growth, and 
replacing $\Rhat(f,N)$ with the corresponding Borel sum replaces the
asymptotic relation \eqref{eq.RfN2} with an exact identity.

In a different direction, suppose that

\begin{itemize}
\item[]
$f$ is continuous in the vertical strip $0 \leq \Re(x) \leq 1$, 
holomorphic in its interior, and 
$f(x) =o(e^{2 \pi | \Im(x)|})$ as $|\Im(x)| \to \infty$ in the strip,
uniformly with respect to $\Re(x)$. 
\end{itemize} 
Then, the {\em Abel-Plana} formula states
that (see \cite[Sec.8.3]{O}):

\begin{equation}
\lbl{eq.abelplana}
\sum_{k=1}^N f\left(\frac{k}{N}\right)=N \int_0^1 f(u)du + \frac{1}{2}
( f(1)-f(0)) + i \int_0^\infty \frac{f\left(\frac{iy}{N}\right)-
f\left(1+\frac{iy}{N}\right) - f\left(-\frac{iy}{N}\right)
+ f\left(1-\frac{iy}{N}\right)}{e^{2 \pi y}-1} dy.
\end{equation}

\subsection{What is a resurgent function?}
\lbl{eq.resurgent}

The notion of a resurgent function was 
introduced and studied by \'Ecalle; see \cite{Ec1}. For our purposes,
a resurgent function is one that admits {\em endless analytic continuation}
(expect at a countable set of non-accumulating singular points)
in the complex plane, and is {\em exponentially bounded}, that is,
satisfies an estimate:
\begin{equation}
\lbl{eq.expgrowth}
|f(z)| < C e^{a|z|}
\end{equation}
for large $z$. Examples of resurgent functions are meromorphic 
functions, algebraic functions, or Borel transforms of solutions of 
generic differential equations with analytic coefficients. 
The $n$th coefficient of the Taylor series of a resurgent function
around a regular point has a manifest asymptotic expansion with respect
to $1/n$ that include small
exponential corrections; see for example \cite[Sec.7]{CG1}. This property
of resurgent functions is key in applications to quantum topology, where
a main problem is to show the existence of asymptotic expansions.
For example, an asymptotic expansion of the coefficients of the
power series \eqref{eq.RfN3} is almost trivial (for a fixed function $f$).
On the other hand, the existence of asymptotic expansion for the 
coefficients of $F_{3_1}(x)$ and $F_{4_1}(x)$ (or more generally, 
$F_{a,b,c}(e^{1/x})$) is a highly non-trivial fact that follows from
the resurgence of the Borel transform of $F_{a,b,c}(e^{1/x})$; see \cite{CG2}.

For an introduction to resurgent functions and their properties, we refer
the reader to the survey articles \cite{CNP1,CNP2,D,DP,Ma,Sa} and for a
thorough study, the reader may consult \'Ecalle's original
work \cite{Ec1,Ec2}.
Let us point out, however, that our main results (Theorems \ref{thm.1} 
and  \ref{thm.3} below) do not require any substantial knowledge of 
resurgence.

\subsection{Statement of the results}
\lbl{sub.results}

Let us recall a useful definition.

\begin{definition}
\lbl{def.borelt}
The {\em (formal) Borel transform} of a formal
power series in $1/N$ is a formal power series in $p$ defined by:
\begin{equation}
\lbl{eq.borelt}
\calB : \BC[[N^{-1}]]  \longto \BC[[p]], 
\qquad
\calB\left(\sum_{n=0}^\infty a_n \frac{1}{N^{n+1}} \right)
=\sum_{n=0}^\infty \frac{a_n}{n!} p^n.
\end{equation}
\end{definition}
Let $G_f(p)$ denote the Borel transform of the power series $\Rhat(f,N)$.

\begin{theorem}
\lbl{thm.1}
If $f(x)$ is resurgent and $f'(x)$ is continuous at $x=0,1$ then 
$G_f(p)$ is given by:
\begin{equation}
\lbl{eq.borelG}
G_f(p)= \frac{1}{4 \pi^2} \sum_{n=1}^\infty
 \frac{1}{n^2}  \left( f'\left(1+\frac{p}{2 \pi i n}\right)
+ f'\left(1-\frac{p}{2 \pi i n}\right)
-f'\left(\frac{p}{2 \pi i n}\right)
-f'\left(-\frac{p}{2 \pi i n}\right)  \right)
\end{equation}
In particular, $G_f(p)$ is resurgent with singularities given by 
\begin{equation}
\lbl{eq.singG}
\calN=\{2 \pi i n \w, \,\, 2 \pi i n (\w-1) \,\, | \,\, n \in \BZ^*, 
\w=\text{singularity of} \,\, f'\}.
\end{equation}
\end{theorem}

Let us consider a function $f(x)$ that satisfies the following:

\begin{itemize}
\item[(A1)]
$f$ is resurgent with no 
singularities in the vertical strip $0 \leq \Re(x) \leq 1$, 
and $f(u)=o(e^{2 \pi |\Im(u)|})$ as $|\Im(u)| \to \infty$ in the strip, 
uniformly with respect to $\Re(u)$.
\end{itemize}

Then, we have the following exact form of the Euler-MacLaurin summation
formula.

\begin{proposition}
\lbl{prop.2}
Under the hypothesis (A1), for every $N \in \BN$ we have:
\begin{equation}
\lbl{eq.em3}
\sum_{k=1}^N f\left(\frac{k}{N}\right)=
N \int_0^1 f(s) ds + \frac{1}{2}(f(1)-f(0)) + \int_0^\infty e^{-Np} G_f(p) dp.
\end{equation}
In particular, the left-hand side of the above equation is the evaluation
at $N$ of an analytic function in the right hand plane.
\end{proposition}

Our proof of Proposition \ref{prop.2} allows to generalize to the case that
$f$ is resurgent with singularities in the vertical strip 
$0 \leq \Re(x) \leq 1$; see Theorem \ref{thm.3} in Section \ref{sub.thm3}.
In that case, every singularity $\l$ of $f$ in the vertical strip
gives rise to exponentially small corrections, and the right hand side of
Equation \eqref{eq.em3} is replaced by a {\em transseries}.

Finally, let us give an integral formula for $G_f(p)$ which is useful
in studying the behavior of $G_f(p)$ for large $p$.

\begin{theorem}
\lbl{thm.11}
With the assumptions of Theorem \ref{thm.1} we have:
\begin{equation}
\lbl{eq.2int}
\displaystyle 
G_f(p)=\frac{1}{(2 \pi i )^3} 
\int_0^\infty \int_{\ga_0} \frac{u}{e^u-1} \left(
\frac{f(s)}{s^2} (e^{\frac{pu}{2 \pi i s}}+e^{-\frac{pu}{2 \pi i s}}) -
\frac{f(1+s)}{(1+s)^2} (e^{\frac{pu}{2 \pi i (1+s)}}+e^{-\frac{pu}{2 \pi i (1+s)}})
\right)
ds du
\end{equation}
where $\ga_0$ is a small circle around $0$ oriented
counterclockwise.
\end{theorem}
Let us end the introduction with some remarks.

\begin{remark}
\lbl{rem.newconstruction}
Theorem \ref{thm.1}, and especially Theorem \ref{thm.3} below
provide a new construction of resurgent functions.
Best known resurgent functions are those that satisfy a difference or
differential equation, linear or not; see for example \cite{Br, BrK,C1} 
and \cite{Ec2}.

 On the other hand, due to the position and shape of their
singularities, the resurgent functions $G_f(p)$ of 
Theorem \ref{thm.1} do not seem to satisfy any differential equations 
with polynomial coefficients.

For example, consider the function $f(x)=(x-\w)^{-m}$ where $\w \not\in [0,1]$
which satisfies the linear differential equation with polynomial coefficients:
$$
(x-\w) f'(x)-mf(x)=0.
$$
$f$ is resurgent, with only one singularity at $x=\w$.
The corresponding resurgent function
$G_f(p)$ of Theorem \ref{thm.1} $G_f(p)$ has infinitely many
singularities on the rays $2 \pi i \w \BR^+$,  $2 \pi i \w \BR^-$,
$2 \pi i (\w-1) \BR^+$, $2 \pi i (\w-1) \BR^-$. It seems unlikely that
$G_f(p)$ satisfies a linear (or a nonlinear) differential equation with 
polynomial coefficients.
\end{remark}

\begin{remark}
\lbl{rem.sing}
Let us point out that Theorem \ref{thm.1} implies 
that the shape of the singularities of $G_f(p)$ is the same as that of $f'(x)$.
For example, if $f'(x)$ is simply ramified,  then so is $G_f(p)$.
We recall that a resurgent function $h(x)$ is {\em simply ramified}
if, locally, at each singularity $\w$ of $h(x)$ we have:
\begin{equation}
\lbl{eq.simplesing}
h(x)=P\left(\frac{1}{x-\w}\right) + \frac{1}{2 \pi i} \log(x-\w) r(x-\w)
+ s(x-\w)
\end{equation}
where $s,r$ are convergent germs, and $P$ is a polynomial.
\end{remark}

\subsection{Acknowledgement}
An early version of this paper was presented by the second author
at talks in Columbia University, Universit\'e Paris VII and Orsay
in the spring and fall of 2006. 
The authors wish to thank J. \'Ecalle for encouraging conversations.
The second author wishes to thank G. Masbaum, W. Neumann, D. Thurston
for their hospitality.

\section{Proof of Theorem \ref{thm.1}}
\lbl{sec.thm1}

\subsection{Computation of the Borel transform $G_f(p)$}
\lbl{sub.computationborel}

Let $\circledast$ denote the {\em Hadamard product} of power series:
\begin{equation}
\lbl{eq.hadamard}
\left( \sum_{n=0}^\infty b_n p^n \right) \circledast
\left( \sum_{n=0}^\infty c_n p^n \right)
=
 \sum_{n=0}^\infty b_n c_n p^n .
\end{equation}
It is classical, and easy to check, 
that the Hadamard product $A \circledast B$
of two functions $A(p)$ and $B(p)$ analytic at $p=0$ is also given by
an integral formula:

\begin{equation}
\lbl{eq.hadamardint}
(A \circledast B)(p)=\frac{1}{2 \pi i} \int_{\gamma} A(s) 
B\left(\frac{p}{s}\right) ds,
\end{equation}
where $\gamma$ is a suitable contour around the origin.
For a detailed explanation of the above formula, see \cite[p.302]{Ju}
and also \cite[p.245]{Bo}.

Let $G_f(p)$ denote the formal Borel transform of the power series in 
\eqref{eq.RfN2}. Since $B_m=0$ for odd $m > 1$, we have:

\begin{eqnarray*}
G_f(p) &=& \calB \left(
\sum_{n=1}^\infty \frac{B_{2n}}{(2n)!} \left( f^{(2n-1)}(1)-
f^{(2n-1)}(0) \right) \frac{1}{N^{2n-1}} \right) \\
&=&
\sum_{n=1}^\infty \frac{B_{2n}}{(2n)!} \left( f^{(2n-1)}(1)-
f^{(2n-1)}(0) \right) \frac{p^{2n-2}}{(2n-2)!} \\
&=&
\sum_{m=2}^\infty \frac{B_{m}}{m!} \left( f^{(m-1)}(1)-
f^{(m-1)}(0) \right) \frac{p^{m-2}}{(m-2)!} \\
&=&
\left( \sum_{m=2}^\infty \frac{B_{m}}{m!} p^{m-2} \right)
\circledast
\left( \sum_{m=2}^\infty \left( f^{(m-1)}(1)-
f^{(m-1)}(0) \right) \frac{p^{m-2}}{(m-2)!} \right) \\
&=& g_1(p) \circledast g_2(p)
\end{eqnarray*}
where
\begin{eqnarray}
\lbl{eq.g1ber}
g_1(p) &=& \sum_{m=2}^\infty \frac{B_{m}}{m!} p^{m-2} =
\frac{1}{p}\left( \frac{1}{e^p-1} -\frac{1}{p} + \frac{1}{2} \right)
\end{eqnarray}
and
\begin{eqnarray}
\lbl{eq.g2p}
g_2(p) &=&
\sum_{m=2}^\infty \left( f^{(m-1)}(1)-
f^{(m-1)}(0) \right) \frac{p^{m-2}}{(m-2)!} \\ &=& \notag
\frac{d}{dp} \left( \sum_{m=2}^\infty \left( f^{(m-1)}(1)-
f^{(m-1)}(0) \right) \frac{p^{m-1}}{(m-1)!} \right) \\
\notag
&=&  \frac{d}{dp}  \left(f(1+p)-f(p)-f(1)+f(0) \right) 
= f'(1+p)-f'(p).
\end{eqnarray} 
Consider positive numbers $r_0$ and $\d$ such that $g_1(p)$ is 
analytic for $|p| < r_0$ (eg, $r_0 < 2 \pi$) and
$g_2(p)$ is analytic for $|p| < \d$--the latter is possible by 
Equation \eqref{eq.g2p} and our assumptions on $f$.

Now, 
Equation \eqref{eq.hadamardint} implies that

\begin{equation}
\lbl{eq.cauchy}
G_f(p)=\frac{1}{2 \pi i} \int_{\gamma} g_1(s) g_2\left( \frac{p}{s} \right)
\frac{ds}{s}
\end{equation}
for all $p$ with $|p| < \d r$,
where $\gamma$ is a circle around $0$ with radius $r<r_0$.

With our assumptions, when $|p| < \d r$, 
the function $s \mapsto g_2(p/s)$ has no singularities outside
of $\gamma$. 
Thus, outside of $\gamma$, the singularities of $g_1(s) g_2(p/s)/s$
are simple poles at the points $2 \pi i n$ for $n \in \BZ^*$, 
with residues
\begin{eqnarray*}
\Res \left( g_1(s) g_2\left(\frac{p}{s} \right)\frac{1}{s}, 
s=2 \pi i n \right) &=& 
\Res(g_1(s), s=2 \pi i n) \,\, g_2\left(\frac{p}{2 \pi i n} \right)
\frac{1}{2 \pi i n} \\
&=& -\frac{1}{4 \pi^2 n^2}  g_2\left(\frac{p}{2 \pi i n} \right).
\end{eqnarray*}
Moreover, $g_1(s)=O(1/s)$ when the distance of $s$ from $2 \pi i \BZ^*$ is
greater than $0.1$ and $g_2(p/s)=O(1)$ for $s$ large, thus
the integrand vanishes at infinity.

We now enlarge the circle $\gamma$ and collect the corresponding
residues by Cauchy's theorem. 
Using the above calculation of the residue and Equation
\eqref{eq.g2p}, it follows that

\begin{eqnarray*}
G_f(p)&=& \frac{1}{4 \pi^2} \sum_{n=1}^\infty
 \frac{1}{n^2}  \left( g_2\left(\frac{p}{2 \pi i n} \right)+
g_2\left(-\frac{p}{2 \pi i n} \right) \right) \\
&=&
\frac{1}{4 \pi^2} \sum_{n=1}^\infty
 \frac{1}{n^2}  \left( f'\left(1+\frac{p}{2 \pi i n}\right)
+ f'\left(1-\frac{p}{2 \pi i n}\right)
-f'\left(\frac{p}{2 \pi i n}\right)
-f'\left(-\frac{p}{2 \pi i n}\right)  \right)
\end{eqnarray*}

Since $f'(x)$ is regular at $x=0,1$, it follows that the above series
is convergent for $p \in \BC-\calN$, where $\calN$ is defined in 
\eqref{eq.singG}. In addition, we conclude that $G_f(p)$ has endless
analytic continuation with singularities in $\calN$.

It remains to prove that $G_f(p)$ is exponentially bounded, assuming that
$f$ is. If $f$ is exponentially bounded, Cauchy's formula implies 
that $f'$ is exponentially bounded. Then, we have:
$$
\left|f'\left(1+\frac{p}{2 \pi i n}\right)\right| \leq
C \exp\left(a\left|1+\frac{p}{2 \pi i n}\right|\right) 
\leq C e^a \exp\left( a\frac{|p|}{2 \pi}\right).
$$
Thus,
$$
|G_f(p)| \leq \frac{C (e^a+1) }{2 \pi^2} \exp\left(a\frac{|p|}{2 \pi}\right) 
\sum_{n=1}^\infty \frac{1}{n^2}=
\frac{C (e^a+1) }{12} \exp\left(a\frac{|p|}{2 \pi}\right).
$$
This completes the proof of Theorem \ref{thm.1}.
\qed

\section{An exact form of Euler-Maclaurin summation formula}
\lbl{sec.exact}

\subsection{Proof of Proposition \ref{prop.2}}
\lbl{sub.thm2}

Proposition \ref{prop.2} follows easily from the Abel-Plana formula;
see Appendix \ref{app.A}. However, we give a proof of Proposition 
\ref{prop.2} that allows us to generalize to Theorem \ref{thm.3} below.

Consider a resurgent function $f$ that satisfies the assumptions (A1),
and let us introduce the function 
$$
h(u)=
\frac{N}{2} f(u) 
\frac{e^{\pi i N u}+e^{-\pi i N u}}{e^{\pi i N u}-e^{-\pi i N u}}
$$
and the contour $\Ga_{R,\d}$ which is a rectangle oriented counterclockwise
with vertices $-iR, 1-iR, 1+iR, iR$ that excludes the points $0$, $1$ 
together with small semicircles of radius $\d$ at the points $0$ and $1$. 

\begin{figure}[htpb]
$$
\psdraw{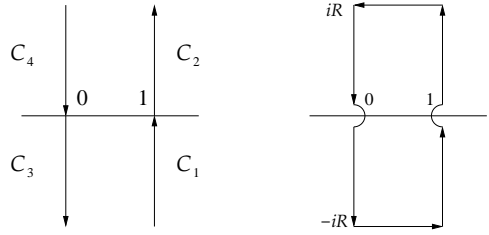}{2.5in}
$$
\caption{The contours $C_1,C_2,C_3,C_4$ of the critical strip, and a 
truncated contour $\Ga_{R,\d}$}\lbl{f.strip}
\end{figure}

Due to our assumptions on $f$, the singularities
of $h(u)$ inside $\Ga_{R,\d}$ are simple
poles at $k/N$ with residue $f(k/N)/(2 \pi i)$ for $k=1,\dots,N-1$.
The residue theorem implies that

\begin{equation}
\lbl{eq.p1}
\sum_{k=1}^n f\left(\frac{k}{N}\right)=\frac{N}{2}
\int_{\Ga_{R,\d}} f(u) 
\frac{e^{\pi i N u}+e^{-\pi i N u}}{e^{\pi i N u}-e^{-\pi i N u}}
du.
\end{equation}
Let $\Ga^+_{R,\d}$ (resp. $\Ga^-_{R,\d}$) denote the upper (resp. lower) part
of the contour $\Ga$. Since $f(x)$ has no singularities in $\Re(u) \in [0,1]$,
the residue theorem implies that

\begin{equation}
\lbl{eq.p2}
-N \int_0^1 f(u)du=\frac{N}{2} \int_{\Ga^+_{R,\d}} f(u)du -\frac{N}{2}
\int_{\Ga^-_{R,\d}} f(u) du.
\end{equation}
Adding up Equations \eqref{eq.p1}, \eqref{eq.p2} and using
\begin{eqnarray*}
\frac{1}{2} \frac{z+z^{-1}}{z-z^{-1}}+\frac{1}{2} &=& \frac{1}{1-z^{-2}} \\
\frac{1}{2} \frac{z+z^{-1}}{z-z^{-1}}-\frac{1}{2} &=& \frac{1}{z^{2}-1} \\
\end{eqnarray*}
we obtain that

\begin{equation}
\lbl{eq.p3}
\sum_{k=1}^{N-1} f\left(\frac{k}{N}\right) -N \int_0^1 f(u)du=
N \int_{\Ga^+_{R,\d}} \frac{f(u)}{1-e^{-2 \pi i N u}} du +
N \int_{\Ga^-_{R,\d}} \frac{f(u)}{e^{2 \pi i N u}-1} du 
\end{equation}

Now let $R \to \infty$. Due to assumption (A1), the integrals over the 
horizontal parts of $\Ga^{\pm}_{R,\d}$ approach zero. 
Next, let $\d \to 0$. Since $f$ is continuous, the integral around
the quarter circle that links $\d$ to $i \d$ tends to $-f(0)/4$.
The other quarter circles are treated similarly.

Thus, we have:

\begin{multline}
\lbl{eq.p4}
\sum_{k=1}^{N} f\left(\frac{k}{N}\right) -N \int_0^1 f(u)du =
\frac{1}{2}(f(1)-f(0)) \\ 
\notag
+
\int_{C_2} \frac{f(u)-f(1)}{1-e^{-2 \pi i N u}} du 
+ \int_{C_4} \frac{f(u)-f(0)}{1-e^{-2 \pi i N u}} du 
+
\int_{C_1} \frac{f(u)-f(1)}{e^{2 \pi i N u}-1} du 
+
\int_{C_3} \frac{f(u)-f(0)}{e^{2 \pi i N u}-1} du. 
\end{multline}

Consider now the corresponding function $G_f(p)$ from Theorem \ref{thm.1}.
We have:
$$
G_f(p)=G_1(p)+G_2(p)+G_3(p)+G_4(p)
$$
where 
\begin{eqnarray*}
G_1(p) &=& \frac{1}{4 \pi^2} \sum_{n=1}^\infty
 \frac{1}{n^2}   f'\left(1+\frac{p}{2 \pi i n}\right) \\
G_2(p) &=& \frac{1}{4 \pi^2} \sum_{n=1}^\infty
 \frac{1}{n^2}   f'\left(1-\frac{p}{2 \pi i n}\right) \\
G_3(p) &=&- \frac{1}{4 \pi^2} \sum_{n=1}^\infty
 \frac{1}{n^2}   f'\left(\frac{p}{2 \pi i n}\right) \\
G_4(p) &=& -\frac{1}{4 \pi^2} \sum_{n=1}^\infty
 \frac{1}{n^2}   f'\left(-\frac{p}{2 \pi i n}\right). 
\end{eqnarray*}
Consider the contours $C_1,C_2,C_3,C_4$ on the boundary of our strip,
as shown in Figure \ref{f.strip}.

We claim that the Laplace transform of the $G_j(p)$ for $j=1,\dots,4$
is given by:

\begin{equation}
\lbl{eq.lapgj}
\int_0^\infty e^{-Np} G_j(p)= 
\begin{cases}
\displaystyle 
N \int_{C_j} \frac{f(u)-f(1)}{e^{2 \pi i N u}-1} du & \,\, \text{for} \,\,
j=1 \\
\displaystyle 
N \int_{C_j} \frac{f(u)-f(0)}{e^{2 \pi i N u}-1} du & \,\, \text{for} \,\,
j=3 \\
\displaystyle 
N \int_{C_j} \frac{f(u)-f(1)}{1-e^{-2 \pi i N u}} du 
& \,\, \text{for} \,\,
j=2 \\ 
\displaystyle 
N \int_{C_j} \frac{f(u)-f(0)}{1-e^{-2 \pi i N u}} du 
& \,\, \text{for} \,\,
j=4.
\end{cases}
\end{equation}
Let us show this for $j=3$; the other integrals are treated in the same
way. We compute as follows:

\begin{xalignat*}{2} 
\int_0^\infty e^{-Np} G_3(p) &=
- \frac{1}{4 \pi^2} \int_0^\infty \sum_{n=1}^\infty  e^{-Np}
 \frac{1}{n^2}   f'\left(\frac{p}{2 \pi i n}\right) dp & \text{by
interchanging sum and integral} \\
&=
 \frac{1}{2 \pi i} \int_{C_3} \sum_{n=1}^\infty 
\frac{e^{-2 \pi i N n u}}{n}   f'(u) du & \text{by} \, p=2 \pi i n u \\
&=
 -\frac{1}{2 \pi i} \int_{C_3} \log(1-e^{-2 \pi i N  u})   
f'(u) du & \text{by \eqref{eq.logsum}} \\
&=
N \int_{C_3}  \frac{f(u)-f(0)}{e^{2 \pi i N  u}-1}
 du & \text{by integration by parts}
\end{xalignat*}
 where
\begin{equation}
\lbl{eq.logsum}
\sum_{n=1}^\infty 
\frac{e^{-2 \pi i N n u}}{n} =-\log(1-e^{-2 \pi i N  u})
\end{equation}
This concludes the proof of Proposition \ref{prop.2} in case $f$ satisfies
(A1).
\qed

Let us end this section with a remark.

\begin{remark}
\lbl{ex.1}
If $f(x)=e^{c x}$, one may verify Equation \eqref{eq.borelG} 
directly by using 
the Mittag-Leffler decomposition of the function $x/(e^x-1)$.
\end{remark}

\subsection{Euler-MacLaurin summation for functions 
with singularities in the vertical strip}
\lbl{sub.thm3}

In this section we consider a function $f$ that satisfies the following
assumptions:
\begin{itemize}
\item[(A2)]
$f$ is resurgent, and let $\La$ denote its set of  singularities
on the critical strip $0 \leq \Re(x) \leq 1$. We assume that $\l \not\in
[0,1]$ for all $\l \in \La$, and $f(u)=o(e^{2 \pi |\Im(u)|})$ as 
$|\Im(u)| \to \infty$ in the strip, uniformly with respect to $\Re(u)$.
We also assume that on a vertical ray $\l+i\BR^+$, we have
$f(u)e^{-2 \pi \Im(u)} \in L^1(\l+i\BR^+)$. 
\item[(A3)]
For every $\l \in \La$ there exist a holomorphic germ $h_{\l}(u)$
and real numbers $\a_{\l}, \b_{\l}$ so that for $u$ near $0$ we have:
\begin{equation}
\lbl{eq.algsing}
f(u+\l)=u^{\a_{\l}} (\log u)^{\b_{\l}} h_{\l}(u).
\end{equation}
\item[(A4)]
For simplicity, let us also assume that
$\Re(\l) \neq \Re(\l')$ for $\l \neq \l'$, and that $\La$ is a finite set.
\end{itemize}

Let 
\begin{equation}
\lbl{eq.laplace}
(\calL G)(x)=\int_0^x e^{-xp} G(p)dp
\end{equation}
denote the {\em Laplace transform} of $G(p)$.
We denote by $f_{\l}(u)$ the {\em variation} (or jump) of the multivalued 
function $f(u+\l)$ at $u$; where $u$ lies the vertical ray starting at $0$
(see for example, \cite{Ma}). We also define: 
\begin{equation}
\lbl{eq.Gsm}
G_{f,\l,m}(p)=i \frac{1}{2 \pi m} f_{\l}\left(\frac{i p}{2 \pi m}\right).
\end{equation}
In case $f(u+\l)$ is single-valued then $G_{f,\l,m}(p)$ is a distribution 
supported at $p=0$.

Then, we have the following exact form of the Euler-MacLaurin summation
formula.

\begin{theorem}
\lbl{thm.3}
\rm{(a)} If $f$ satisfies (A1-A4) and $\a_{\l} > -1$ for all $\l \in \La$, 
then for every $N \in \BN$ we have:
\begin{eqnarray}
\lbl{eq.em4}
\sum_{k=1}^N f\left(\frac{k}{N}\right) &=&
N \int_0^1 f(s) ds + \frac{1}{2}(f(1)-f(0)) + 
(\calL G_f)(N) \\
\notag
& + & 
 N e^{2 \pi i \l  N} \sum_{\l: \Im(\l)>0} 
\sum_{m=0}^\infty e^{2 \pi i \l m N} (\calL G_{f,\l,m})(N) \\
\notag
& + & N e^{-2 \pi i \l N}
\sum_{\l: \Im(\l)<0} 
\sum_{m=0}^\infty e^{-2 \pi i \l m N} (\calL G_{f,\l,m})(N) 
\end{eqnarray}
\rm{(b)} If some $\a_{\l} \leq -1$, Equation \eqref{eq.em4} is true
after integration by parts $M$-times where 
$M \geq \max_{\l: \a_{\l} \leq -1}[-\a_{\l}]$.
\end{theorem}

\begin{proof}
Without loss of generality, let us assume that $f$ has a single singularity
$\l$ in the vertical strip $0 \leq \Re(x) \leq 1$ with $\Im(\l)>0$.

Use the modified contour $\Ga_{R,\d,\l}$ in Figure \ref{f.strip2}.
\begin{figure}[htpb]
$$
\psdraw{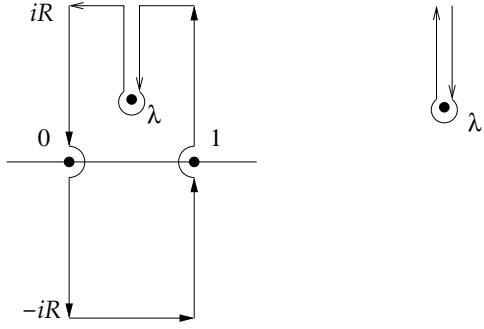}{2.5in}
$$
\caption{The modified truncated contour $\Ga_{R,\d,\l}$ on the left and
a Hankel contour $H$ on the right.}\lbl{f.strip2}
\end{figure}

Let $H_R$ denote the portion of $\Ga_{R,\d,\l}$ that consists of the
truncated Hankel contour around $\l$, and $S_{R,\d}=\Ga_{R,\d,\l}-H_R$.
Equations \eqref{eq.p1} and \eqref{eq.p2} become:
\begin{equation}
\lbl{eq.p11}
-\frac{N}{2} \int_{H_R} f(u) 
\frac{e^{\pi i N u}+e^{-\pi i N u}}{e^{\pi i N u}-e^{-\pi i N u}}
du +
\sum_{k=1}^n f\left(\frac{k}{N}\right)=\frac{N}{2}
\int_{S_{R,\d}} f(u)
\frac{e^{\pi i N u}+e^{-\pi i N u}}{e^{\pi i N u}-e^{-\pi i N u}}
du
\end{equation}
and

\begin{equation}
\lbl{eq.p22}
-\frac{N}{2} \int_{H_R} f(u)du
-N \int_0^1 f(u)du=\frac{N}{2} \int_{S^+_{R,\d}} f(u)du -\frac{N}{2}
\int_{S^-_{R,\d}} f(u) du.
\end{equation}
Adding up, the extra contribution from $H_R$ becomes:

\begin{equation}
\lbl{eq.p33}
-N \int_{H_R} \frac{f(u)}{1-e^{-2 \pi i N u}} du
\end{equation}
Now let $R \to \infty$. 
Notice that $f(u+\l)$ is uniformly $L^1$ for $u$ near $0$ iff $\a_{\l} > -1$
for all $\l$.  Using this and our integrability assumption (A2),
it follows that in the 
limit the above integral equals to

\begin{eqnarray*}
I_\l&=& -N \int_{H} \frac{f(u)}{1-e^{-2 \pi i N u}} du \\
&=&
-N i \int_0^\infty \frac{f_\l(is)}{1-e^{-2 \pi i N (\l+is)}}  ds.
\end{eqnarray*}
Now, $\l=\l_1+i \l_2$ with $\l_2 >0$, and we may write
\begin{eqnarray*}
\frac{1}{1-e^{-2 \pi i N (\l+is)}} &=& \frac{1}{1-\omega^{-N} e^{2 \pi N s}} \\
& =& - \omega^N e^{-2 \pi N s} \sum_{m=0}^\infty \omega^{N m} e^{-2  N m s}
\end{eqnarray*}
where $\omega=e^{2 \pi i \l}$ satisfies $|\omega| <1$.
Thus, 

\begin{eqnarray*}
I_\l=-N \omega^{N} \sum_{m=0}^\infty \omega^{N m} (\calL G_{\l,m})(N). 
\end{eqnarray*}
Part (a) of Theorem \ref{thm.3} follows. Part (b) follows from the fact that
if $f(u+\l)$ has a local expansion of the form \eqref{eq.algsing}, and
$F^{(s)}(u)=f(u)$, then $F(u+\l)$ as a local expansion of the form:

\begin{equation}
\lbl{eq.algsing2}
F(u+\l)=u^{\a_{\l}+s} (\log u)^{\b_{\l}} H_{\l}(u)
\end{equation}
for a holomorphic germ $H_{\l}(u)$. Cf. also \cite[Thm.1]{C1}.
\end{proof}

\subsection{Euler-MacLaurin with logarithmic singularities at $x=0$}
\lbl{sub.logsing}

In this section we consider functions $f(x)$ that have a logarithmic
singularity at $x=0$. 
Motivated by our applications to quantum topology, we 
consider functions $f$ of the form: 

\begin{equation}
\lbl{eq.flog}
f(x)=c \log x + g(x)
\end{equation}
where $g$ that satisfies (A1), and $c \in \BC$. Let us define

\begin{equation}
\lbl{eq.H}
H(p) = \frac{1}{p^2} \left( \frac{p}{e^p-1}-1+\frac{p}{2} \right)
\end{equation}
It is easy to see that $H(p)$ is analytic at $p=0$. In fact,
the Taylor series of $H$ at $p=0$ is given by:

\begin{equation}
\lbl{eq.stir}
\sum_{n=1}^\infty \frac{B_{2n}}{(2n)!} p^{2n-2}
\end{equation}

\begin{theorem}
\lbl{thm.4}
Under the above hypothesis, for every $N \in \BN$ we have:
\begin{equation}
\lbl{eq.em5}
\sum_{k=1}^N f\left(\frac{k}{N}\right)=
N \int_0^1 f(s) ds + \frac{c}{2} \log N +\frac{c}{2} \log(2 \pi) 
+ \frac{1}{2}(g(1)-g(0)) + \calL(G_g+ c H)(N).
\end{equation} 
\end{theorem}

\begin{proof}
Since $f$ is given by \eqref{eq.flog}, we have:
$$
\sum_{k=1}^N f\left(\frac{k}{N}\right)=
c \log\left(\frac{N!}{N^N}\right)
+ \sum_{k=1}^N g\left(\frac{k}{N}\right)
$$
Recall now from \cite[Sec.13.15]{Ha} the
following exact form of {\em Stirling's formula}:

\begin{equation}
\lbl{eq.stirling}
\log\left(\frac{N!}{N^N}\right) =
\frac{1}{2} \log N -N  +\frac{1}{2} \log(2 \pi) + (\calL H)(N).
\end{equation}
Applying Proposition \ref{prop.2} to $g$ gives:
$$
\sum_{k=1}^N g\left(\frac{k}{N}\right)=
N \int_0^1 g(s) ds + \frac{1}{2}(g(1)-g(0)) + (\calL G_g)(N).
$$
Adding up, and using 
$$
N \int_0^1  f(s) ds= N \int_0^1  g(s) ds + N c \int_0^1  \log s ds 
=N \int_0^1  g(s) ds -N c
$$
we obtain \eqref{eq.em5}. The result follows.
\end{proof}

\section{Parametric resurgence of difference equations with a parameter}
\lbl{sec.parametric}

Consider the first order linear difference equation 
with a small parameter $\e$:
\begin{equation}
\lbl{eq.diff}
y(x+\e,\e)=a(x,\e) y(x,\e)
\end{equation}
where $a(x,\e)$ is smooth. 
\eqref{eq.diff} has a unique {\em formal solution} (often called a {\em
WKB solution}) of the form:

\begin{equation}
\lbl{eq.WKB}
y(x,\e)=e^{\frac{1}{\e} \sum_{k=0}^\infty F_k(x) \e^k }
\end{equation}
where $F_j(0)=0$. See for example, \cite{CC} and \cite{GG}.
For simplicity, suppose that $a(x,\e)=a(x)$ is independent 
of $\e$.
Under the stated assumptions, the next theorem gives an exact solution to 
\eqref{eq.diff} which is asymptotic to the formal solution \eqref{eq.WKB}.

\begin{theorem}
\lbl{thm.diff}
\rm{(a)} For all $x$ such that $s \to \log a(s x)$ satisfies (A1)
we have:
\begin{equation}
\lbl{eq.actual2}
\frac{1}{\e} \sum_{k=0}^\infty F_k(x) \e^k \sim
\frac{1}{\e} \int_0^x \log a(q) dq -\frac{1}{2} \log a(x) + \frac{1}{2}
\log a(0) + \int_0^\infty e^{-q/\e} G(q,x) dq
\end{equation}
where
\begin{multline}
\lbl{eq.borelGG}
G(q,x) = \frac{1}{4 \pi^2} \sum_{n=1}^\infty
 \frac{1}{n^2}  \left( 
\frac{a'\left(x+\frac{q}{2 \pi i n}\right)}{
a\left(x+\frac{q}{2 \pi i n}\right)}
+
\frac{a'\left(x-\frac{q}{2 \pi i n}\right)}{
a\left(x-\frac{q}{2 \pi i n}\right)}
-
\frac{a'\left(\frac{q}{2 \pi i n}\right)}{
a\left(\frac{q}{2 \pi i n}\right)}
-
\frac{a'\left(-\frac{q}{2 \pi i n}\right)}{
a\left(-\frac{q}{2 \pi i n}\right)}
\right) \\
\displaystyle
= 
\frac{1}{(2 \pi i )^3} 
\int_0^\infty \int_{\ga_0} 
\frac{u}{e^u-1} \left(
\frac{\log a(s)}{s^2} (e^{\frac{pu}{2 \pi i s}}+e^{-\frac{pu}{2 \pi i s}})
- 
\frac{\log a(x+s)}{(x+s)^2} (e^{\frac{pu}{2 \pi i (x+s)}}+e^{-\frac{pu}{2 \pi i (x+s)}})
\right)
ds du. 
\end{multline}
where $\ga_0$ is a small circle around $0$ oriented counterclockwise.
\newline
\rm{(b)} Moreover, \eqref{eq.diff} has a solution $y(x,\e)$ of the form:
\begin{equation}
\lbl{eq.actualsol}
y(x,\e)=\sqrt{\frac{a(0)}{a(x)}} \exp\left(
\frac{1}{\e} \int_0^x \log a(q) dq + 
\int_0^\infty e^{-q/\e} G(q,x) dq \right).
\end{equation}
\end{theorem}

\begin{remark}
\lbl{rem.parametricresurgence}
It follows that the singularities of $G(q,x)$ are of the form 
$2 \pi i n \l$
or $2 \pi i n (\l-x)$ where $n \in \BZ^*$ and $\l$ is a singularity of 
$\log a$. These type of singularities appear in parametric 
(i.e., co-equational) resurgence of \'Ecalle; see \cite{Ec2}.

The proof of Theorem \ref{thm.diff} indicates the close relation 
between the Euler-MacLaurin 
summation formula and the formal solutions of a linear difference equation
with a parameter.

From that point of view, resurgence of $G_f(p)$ translates to parametric
resurgence of formal solutions of linear difference equations. In the case
of formal solutions of linear differential equations with a parameter,
\'Ecalle shows that their singularities are of the form $n (\a_i-x)$ for 
$n=-1,1,2,3,\dots$; see \cite[Eqn.(6.9)]{Ec2}. 
\end{remark}

\begin{proof}
(a) Let $z(x,\e)=\log y(x,\e)$.
Taking the logarithm of \eqref{eq.diff}, it follows that
$$
z(k\e+\e,\e)=\log a(k\e) + z(k\e,\e).
$$ 
Summing up for $k=0,\dots, N-1$ and using the variable 
\begin{equation}
\lbl{eq.x=Ne}
x=N \e
\end{equation}
we obtain that:

\begin{eqnarray*}
z(x,\e)-z(0,\e) &=& \sum_{k=0}^{N-1} \log a(k \e) \\
& = &  
 -\log a(x) + \log a(0) + \sum_{k=1}^{N} \log a(x k/N) .
\end{eqnarray*}
Let us fix $x$ and apply Proposition \ref{prop.2} to the function 
$s \to \log a(sx)$. We obtain that

\begin{eqnarray*}
z(x,\e)-z(0,\e) 
&=&  N \int_0^1 \log a(x s)ds -\frac{1}{2} \log a(x) + \frac{1}{2} \log a(0) +
\int_0^\infty e^{-N p} H(p,x) dp \\
&=&  \frac{1}{\e} \int_0^x \log a(s)ds -\frac{1}{2} \log a(x) 
+ \frac{1}{2} \log a(0) + \int_0^\infty e^{-x p/\e} H(p,x) dp \\
&=&  \frac{1}{\e} \int_0^x \log a(s)ds -\frac{1}{2} \log a(x) 
+ \frac{1}{2} \log a(0) + \int_0^\infty e^{-q/\e} H\left(\frac{q}{x},x\right) 
\frac{dq}{x} 
\end{eqnarray*}
where by Theorems \ref{thm.1} and \ref{thm.11} we have:

\begin{multline*}
H(p,x) = \frac{x}{4 \pi^2} \sum_{n=1}^\infty
 \frac{1}{n^2}  \left( 
\frac{a'\left(x\left(1+\frac{p}{2 \pi i n}\right)\right)}{
a\left(x\left(1+\frac{p}{2 \pi i n}\right)\right)}
+
\frac{a'\left(x\left(1-\frac{p}{2 \pi i n}\right)\right)}{
a\left(x\left(1-\frac{p}{2 \pi i n}\right)\right)}
-
\frac{a'\left(x \frac{p}{2 \pi i n}\right)}{
a\left(x\frac{q}{2 \pi i n}\right)}
-
\frac{a'\left(-x\frac{p}{2 \pi i n}\right)}{
a\left(-x\frac{p}{2 \pi i n}\right)}
\right) \\ \displaystyle
=\frac{1}{(2 \pi i )^3} 
\int_0^\infty \int_{\ga_0} 
\frac{u}{e^u-1} \left(
\frac{\log a(sx)}{s^2} (e^{\frac{pu}{2 \pi i s}}+e^{-\frac{pu}{2 \pi i s}}) -
\frac{\log a((1+s)x)}{(1+s)^2} (e^{\frac{pu}{2 \pi i (1+s)}}
+e^{-\frac{pu}{2 \pi i (1+s)}})
\right)
ds du.
\end{multline*}

\noindent
Since $H(q/x,x)/x=G(q,x)$ (where $G(q,x)$ is given by \eqref{eq.borelGG})
and $z(0,\e)=0$, it follows 
that for all $\e>0$ and $k \in \BZ$ we have:
\begin{equation}
\lbl{eq.zke}
z(k\e+\e,\e)=\log a(k\e)+z(k\e,\e).
\end{equation}
To prove (b), let us consider the difference
$$
E(x,\e)=\e( \log y(x+\e,\e)-\log a(x) - \log y(x,\e)).
$$
It follows by definition that $E(x,\e)$ is analytic in $(x,\e)$.
Thus,
$$
E(x,\e)=\sum_{i,j=0}^\infty c_{ij} x^i \e^j.
$$
Moreover, \eqref{eq.zke} implies 
that for all $\e>0$ and all $k \in \BZ$ we have:
$$
0=F(k\e,\e)=\sum_{i,j=0}^\infty c_{ij} k^i \e^{i+j}.
$$
Thus $c_{i,j}=0$ for all $i,j$ and $E(x,\e)=0$. This completes the proof of 
(b).

The definition of $y(x,\e)$ by a Laplace integral and Watson's lemma
(see \cite[Sec.4.3.1]{O}) implies that
$$
y(x,\e) \sim \frac{1}{\e} \sum_{k=0}^\infty \phi_k(x) \e^k
$$
for analytic functions $\phi_k(x)$ that satisfy $\phi_k(0)=0$. Since
a formal WKB solution given by \eqref{eq.WKB} is unique, it follows that
$F_k(x)=\phi_k(x)$ for all $k$. Thus, (a) follows.
\end{proof}

\begin{remark}
\lbl{rem.1}
Theorem \ref{thm.diff} can be generalized when 
$$
a(x,\e)=\sum_{k=0}^\infty a_k(x) \e^k
$$
is analytic with respect to $(x,\e)$, and the coefficients $a_k(x)$
are resurgent functions. It may also be generalized to the case of
higher order linear difference equations with a parameter. 
This will be explained elsewhere.
\end{remark}

\begin{remark}
\lbl{rem.sauzin}
The reader may compare Theorem \ref{thm.diff} with the results of the last
section of \cite{Sa}.
\end{remark}

\section{An integral formula for $G_f(p)$}
\lbl{sec.thm3}

In this section we give a proof of Theorem \ref{thm.11}.
We follow the ideas of \cite{C2} to convert the sum of Equation
\eqref{eq.borelG} into an integral. Let us show that
\begin{equation}
\lbl{eq.1term}
\sum_{n=1}^\infty
 \frac{1}{n^2}  
f'\left(\frac{p}{2 \pi i n}\right)
= \frac{1}{2 \pi i} 
\int_0^\infty \int_{\ga_0} \frac{u f(s)}{s^2 (e^u-1)} e^{\frac{pu}{2 \pi i s}} ds du
\end{equation}
and similarly for the sum of the other three terms in \eqref{eq.borelG}.

To prove Equation \eqref{eq.1term}, we first expand $f'$ at $p=0$, then
take a Laplace transform with respect to the summation variable $n$,
interchange the order of summation and sum the geometric series.
We obtain that:
\begin{eqnarray*}
\sum_{n=1}^\infty
 \frac{1}{n^2}  
f'\left(\frac{p}{2 \pi i n}\right) &=&
\sum_{n=1}^\infty \sum_{j=0}^\infty \frac{f^{(j+1)}(0)}{j!} 
\left(\frac{p}{2 \pi i}\right)^j \frac{1}{n^{j+2}} \\
& = &
\sum_{n=1}^\infty \sum_{j=0}^\infty \frac{f^{(j+1)}(0) p^j}{j!(j+1)!(2 \pi i)^j} 
\int_0^\infty
e^{-nu} u^{j+1} du \\
& = &
\sum_{j=0}^\infty \int_0^\infty \frac{f^{(j+1)}(0) p^j}{j!(j+1)!(2 \pi i)^j} 
 u^{j+1} \frac{1}{e^u-1} du 
\end{eqnarray*}
Using Cauchy's formula
$$
\frac{f^{(j+1)}(0)}{(j+1)!}=
\frac{1}{2 \pi i}\int_{\ga_0} \frac{f(s)}{s^{j+2}} ds
$$
and interchanging summation and integration it follows that:
\begin{eqnarray*}
\sum_{n=1}^\infty
 \frac{1}{n^2}  
f'\left(\frac{p}{2 \pi i n}\right) &=&
\frac{1}{2 \pi i} \sum_{j=0}^\infty \int_0^\infty \int_{\ga_0}
\frac{f(s) p^j}{j! s^{j+2}} 
 \frac{u^{j+1}}{(2 \pi i)^j} \frac{1}{e^u-1} ds du \\
&=& 
\frac{1}{2 \pi i} \int_0^\infty \int_{\ga_0} \frac{u f(s)}{s^2(e^u-1)} 
\sum_{j=0}^\infty \frac{1}{j!} \left(\frac{p u}{2 \pi i s}\right)^j dsdu \\
& = & \frac{1}{2 \pi i} \int_0^\infty \int_{\ga_0} \frac{u f(s)}{s^2(e^u-1)} 
 e^{\frac{pu}{2 \pi i s}} dsdu.
\end{eqnarray*}
The interchanges of summation and integration are justified by dominated
convergence. 
This concludes the proof of \eqref{eq.1term} and Theorem \ref{thm.11}.
\qed

\appendix

\section{}
\lbl{app.A}

For completeness, let us show how the Abel-Plana formula implies 
Proposition \ref{prop.2}.
With the notation as in Proposition \ref{prop.2}, we claim that for every
$N \in \BN$ we have:
\begin{eqnarray}
\lbl{eq.AP1}
-i \int_0^\infty \frac{f\left(1+\frac{iy}{N}\right)-f(1)}{e^{2 \pi y}-1}=
\frac{1}{4 \pi^2} \sum_{n=1}^\infty \frac{1}{n^2} \int_0^\infty e^{-Np} 
f'\left(1-\frac{p}{2 \pi i n}\right) dp \\
\lbl{eq.AP2}
i \int_0^\infty \frac{f\left(1-\frac{iy}{N}\right)-f(1)}{e^{2 \pi y}-1}=
\frac{1}{4 \pi^2} \sum_{n=1}^\infty \frac{1}{n^2} \int_0^\infty e^{-Np} 
f'\left(1+\frac{p}{2 \pi i n}\right) dp \\
\lbl{eq.AP3}
i \int_0^\infty \frac{f\left(\frac{iy}{N}\right)-f(0)}{e^{2 \pi y}-1}=
-\frac{1}{4 \pi^2} \sum_{n=1}^\infty \frac{1}{n^2} \int_0^\infty e^{-Np} 
f'\left(-\frac{p}{2 \pi i n}\right) dp \\
\lbl{eq.AP4}
-i \int_0^\infty \frac{f\left(-\frac{iy}{N}\right)-f(0)}{e^{2 \pi y}-1}=
-\frac{1}{4 \pi^2} \sum_{n=1}^\infty \frac{1}{n^2} \int_0^\infty e^{-Np} 
f'\left(\frac{p}{2 \pi i n}\right) dp.
\end{eqnarray}
Adding up, and using the Abel-Plana formula \eqref{eq.abelplana}, 
gives a proof of Proposition \ref{prop.2}.
Let us give the proof of \eqref{eq.AP1} and leave the rest as an exercise.
For $y>0$, we have $e^{-2 \pi y} < 1$ and the geometric series gives:

\begin{equation}
\lbl{eq.geomseries}
\frac{1}{e^{2 \pi y}-1}=\sum_{n=1}^\infty e^{-2 \pi n y}.
\end{equation} 
Interchanging summation and integration, changing variables $2 \pi n y=Np$
and integrating by parts
(justified by the hypothesis (A1)), we obtain that

\begin{eqnarray*}
-i \int_0^\infty \frac{f\left(1+\frac{iy}{N}\right)-f(1)}{e^{2 \pi y}-1} &=&
-i \sum_{n=1}^\infty \int_0^\infty e^{-2 \pi n y}
\left(f\left(1+\frac{iy}{N}\right)-f(1) \right) dy \\
&=&
-\frac{iN}{2 \pi} \sum_{n=1}^\infty \frac{1}{n} \int_0^\infty
e^{- Np}
\left(f\left(1-\frac{p}{2 \pi i n}\right)-f(1) \right) dp \\
&=&
\frac{i}{2 \pi} \sum_{n=1}^\infty \frac{1}{n} \int_0^\infty
(e^{- Np})'
\left(f\left(1-\frac{p}{2 \pi i n}\right)-f(1) \right) dp \\
&=&
\frac{1}{4 \pi^2} \sum_{n=1}^\infty \frac{1}{n^2} \int_0^\infty
e^{- Np}
f'\left(1-\frac{p}{2 \pi i n}\right) dp. 
\end{eqnarray*}

\ifx\undefined\bysame
    \newcommand{\bysame}{\leavevmode\hbox
to3em{\hrulefill}\,}
\fi

\end{document}